\documentclass{amsart}
\usepackage{amssymb,amsmath,amscd,xypic}
\newtheorem{theorem}{Theorem}[section]
\newtheorem{lemma}[theorem]{Lemma}
\newtheorem{corollary}[theorem]{Corollary}
\newtheorem{definition}[theorem]{Definition}

\newtheorem{remark}[theorem]{\it Remark}

\newtheorem{proposition}[theorem]{Proposition}

\begin{document}

\large

\def\sMr{\widetilde{M}_\mathbf{r}(\mathbb{C}P^m)}
\def\sMs{\widetilde{M}_{\Lambda-\mathbf{r}}(\mathbb{C}P^{n-m-2})}
\def\Mr{M_\mathbf{r}(\mathbb{C}P^m)}
\def\Nr{N_\mathbf{r}(\mathbb{C}P^m)}
\def\Ms{M_{\Lambda - \mathbf{r}}(\mathbb{C}P^{n-m-2})}
\def\Ns{N_{\Lambda - \mathbf{r}}(\mathbb{C}P^{n-m-2})}
\def\Gr{\mathrm{Gr}}
\def\Grk{Gr_k(\mathbb{C}^n)}
\def\b{\mathbb}
\def\C{\mathbb{C}}
\def\CP{\mathbb{CP}}
\def\Q{\mathbb{Q}}
\def\R{\mathbb{R}}
\def\Z{\mathbb{Z}}
\def\N{\mathbb{N}}
\def\br{\mathbf{r}}
\def\bs{\mathbf{s}}
\def\bk{\mathbf{k}}
\def\bl{\mathbf{l}}
\def\t{\mathrm}
\def\ba{\mathbf{a}}
\def\bb{\mathbf{b}}
\def\Lf{L_\bk^\ba}
\def\Lg{L_\bl^\bb}
\def\Fk{F_{\bk}(\C^n)}
\def\Fl{F_{\bl}(\C^n)}
\def\Rk{R_{\bk}}
\def\Rl{R_{\bl}}
\def\L{\mathbf{\Lambda}}
\def\Zo{\Z_{\scriptscriptstyle{\geq 0}}}
\def\Mu{\mathcal{M}_\br(\mathbb{CP}^m)}
\def\Rlm{R_{\lambda,\mu}}
\def\GL{\mathrm{GL}}
\def\SL{\mathrm{SL}}
\def\PGL{\mathrm{PGL}}
\def\state{\mathrm{state}}
\def\weight{\mathrm{weight}}
\def\q{/ \! /}
\def\supp{\mathrm{supp}}

\title[Matroids and Geometric Invariant Theory]{Matroids and Geometric Invariant Theory of torus actions
on flag spaces}

\author{Benjamin J. Howard*}
\thanks{*The author was partially supported by NSF grant DMS-0104006}
\thanks{The author intends to include these results in his doctoral
thesis}
\date{\today}

\begin{abstract}
Let $F \q T$ be a G.I.T. quotient of a flag manifold $F$ by the
natural action of the maximal torus $T$ in $\SL(n,\C)$. The
construction of the quotient space depends upon the choice of a
$T$-linearized line bundle $L$ of $F$.  This note concerns the case where $L =
L_\lambda$ is a very ample homogeneous line bundle determined by a
dominant weight $\lambda$.  

A theorem of Gel'fand, Goresky, MacPherson, and Serganova
about matroids and matroid polytopes is applied to study semistability
of flags relative to a given $T$-linearization of $L_\lambda$. The main theorem of this note 
is the follwing:  if there exists a nonzero $T$-invariant global section 
of $L_\lambda$, then for each semistable flag $x \in F$, there exists a $T$-invariant 
global section $s$ of $L_\lambda$ such that $s(x) \neq 0$.  Hence the global $T$-invariant 
sections of $L_\lambda$ determine a well-defined map from $F \q T$ to projective space, 
provided there is at least one such which is nonzero.

A related result in this note is that the closure of any $T$-orbit in $F$ is
projectively normal for any projective embedding of $F$.  The proof of this fact 
uses essentially the same argument given for the semistability theorem above.
\end{abstract}

\maketitle

\tableofcontents

\section{Introduction}

The geometry (both symplectic and algebraic) of the quotients $F \q T$ have 
been extensively studied in recent years; Allen Knutson 
called them ``weight varieties''\footnote{The term ``weight variety'' actually refers to 
more general quotients;  they are G.I.T. quotients of $G/P$ by a maximal 
torus $T$ in $G$, where $G$ is a  reductive connected complex Lie group and $P$ is 
a parabolic subgroup of $G$ containing $T$.}
in his thesis \cite{Knutson}.  The dependence of the geometry  
of the quotient on the choice of linearization was studied by Yi Hu \cite{Hu} 
and the cohomology of weight varieties 
was computed by Rebecca Goldin \cite{Goldin}.  Special cases of weight varieties have 
been studied since the nineteenth century; for example a G.I.T. quotient 
$(\CP^{k-1})^n \q \PGL(k,\C)$ 
is isomorphic to a G.I.T. quotient $\Gr_k(\C^n) \q T$ by the Gel'fand MacPherson 
correspondence (here $\Gr_k(\C^n)$ denotes the Grassmannian).  The projective invariants 
of $n$-tuples of points on projective space are still not understood today; we do not 
even know a minimal set of generators for the ring of projective invariants (see page 8 of \cite{Harris}).
 
We take one step towards a solution 
to the generators problem (for $G = \SL(n,\C)$) with Theorem \ref{maintheorem}, which implies that 
the lowest degree $T$--invariants in the graded ring of $F$ 
are sufficient to give a well-defined map from $F \q T$ to projective space.
Consequently these global sections determine an ample line bundle $M$ of $F \q T$.  
We are left with the problem of determining which tensor power of $M$ is very ample.

The proof of Theorem \ref{maintheorem} involves a simple combinatorial
argument involving Minkowski sums of weight polytopes of flags.  These 
weight polytopes are also known as flag matroid polytopes, see \cite{BorovikGelfandWhite}.
Two facts are essential to the argument:
\begin{itemize}
\item Any subset of $\SL(n,\C)$ roots which are linearly independent may be extended to 
a basis of the root lattice. \item Each edge of a matroid polytope
is parallel to a root of $\SL(n,\C)$ (due to Gel'fand, Goresky,
MacPherson, Serganova).
\end{itemize}
\begin{remark}
The first fact is specific to $\SL(n,\C)$.  The root systems of
other classical complex simple Lie algebras do not have this
remarkable saturation property.  However, the second result is a
special case of the Gel'fand--Serganova theorem which is one of
the central theorems in the new subject of Coxeter matroids, see
\cite{BorovikGelfandWhite}.  
It should be noted that Theorem \ref{maintheorem} easily follows
from a theorem of Neil White \cite{White} in the case that $F$ is
a Grassmannian.
\end{remark}

Additionally, the tools we develop in proving Theorem \ref{maintheorem} 
also allow us to show that the closure of a $T$-orbit $cl(T
\cdot x)$ for any $x \in F$ (for any projective
embedding of $F$) is a projectively normal toric variety. Again
Neil White \cite{White} showed this holds when $F$ is a
Grassmannian $\Gr_k(\C^n)$. Additionally, R. Dabrowski
\cite{Dabrowski} proved that projective normality holds for
closures of certain \emph{generic} $T$-orbits in other homogeneous
spaces $G/P$ (he covered the case that $G$ is any semi-simple
complex Lie group). 

\vskip 12pt

{\bf Acknowledgements} The author thanks his advisor John J.
Millson for his advice and patience.   Also, he thanks Professors
Joseph Bonin and Edward Swartz for many useful conversations about
matroids.  Additionally he thanks Professors Allen Knutson and 
Shrawan Kumar for suggestions that lead to the final version of this note.

\section{The construction of $F \q T$ and statement of main theorem}

A weight variety of $G = \SL(n,\C)$ is a G.I.T. quotient of a flag
manifold $F = G/P$ by the action of the Cartan subgroup $T$. The
construction of such a quotient involves the choice of a
$T$-linearized line bundle $L$ of $F = G/P$.  If $L$ is very ample,
then it's isomorphism class is determined by a choice of dominant
weight $\lambda$ such that $P$ is the largest parabolic subgroup
such that the character $e^\lambda$ defined on the Borel subgroup
$B$ of upper triangular matrices extends (uniquely) to $P$. The
$T$-linearization of $L$ will also depend on a choice of a weight
$\mu$, but $\mu$ need not be dominant.

\subsection{Elementary notions from the representation theory of $\SL(n,\C)$}

Since $\SL(n,\C)$ is simply connected, the set of $\SL(n,\C)$
weights are the differentials evaluated at the identity element of
characters $\chi : T \to \C^\ast$, which are holomorphic
homomorphisms (that is, the character lattice coincides with the
weight lattice). The differential $d \chi$ (evaluated at the
identity element of $T$) of $\chi$ lies within the dual Lie algebra
$\mathfrak{t}^\ast$ of $T$. On the other hand, if $\varpi \in
\mathfrak{t}^\ast$ is a weight, we shall denote $e^\varpi$ as the
unique character $e^\varpi : T \to \C^\ast$ such that $d(e^\varpi) =
\varpi$.

A character $e^\lambda$ applied to $t =
\text{diag}(t_1,\ldots,t_n) \in T$ must be equal to $\prod_{i=1}^n
t_i^{a_i}$ for some fixed integers $a_i$. Since $\prod_{i=1}^n t_i
= 1$ for all $t \in T$, we have that the $n$-tuple of exponents
$(a_1,\ldots,a_n)$ and $(a_1+a,a_2+a,\ldots,a_n+a)$ determine the
same character.  We may thus view the abelian group of characters
of $T$ as $\Z^n / \Delta$ where $\Delta$ is the diagonal.  On the
other hand, the weight $\lambda \in \mathfrak{t}^\ast$ takes a
complex vector $(z_1,\ldots,z_n) \in \mathfrak{t}$ (where $z_1 +
z_2 + \cdots + z_n = 0$) to $\sum_{i=1}^n a_i z_i$.  Again, adding
a constant to each $a_i$ results in the same function, and so
again we have that the additive group of weights is isomorphic to
$\Z^n/\Delta$.  We shall henceforth identify characters and
weights as $n$-tuple of integers modulo the diagonal $\Delta$.

\subsubsection{Dominant weights and construction of very ample
line bundles}

We say that a weight $\lambda = (a_1,\ldots,a_n)$ is dominant if
the $a_i$'s are non-strictly decreasing.  Now suppose that
$\lambda$ is dominant, and $P \subset G$ is the largest parabolic
subgroup (a subgroup containing all the upper triangular matrices
in $G$) such that $e^\lambda$ extends to a character $\chi : P \to
\C^\ast$.  It is a basic fact that $\chi$ is determined by it's
restriction $e^\lambda$ to the torus $T$, so we will abuse
notation and identify $\chi$ with $e^\lambda$.

The dominant weight $\lambda$ determines a very ample line bundle
$L_{\lambda}$ of $G/P$.  The total space of $L_\lambda$ is the set
of equivalence classes of pairs $(g,z)$ for $g \in G$ and $z \in
\C$, where $(g,z) \sim (gp, e^\lambda(p)z)$ for all $p \in P$. The
map $\pi$ from the total space to $G/P$ is given by $\pi : (g,z)
\mapsto gP$.  Each global section of $L_\lambda$ is given by
$s_f(gP) = (g,f(g))$ where $f : G \to \C$ is a holomorphic function
such that $f(gp) = e^\lambda(p)f(g)$ for all $p \in P$ and $g \in
G$.

There is a natural action of $G$ on the total space of
$L_\lambda$, given by $g \cdot (g',z) = (gg', z)$. This defines an
action on sections by $$(g \cdot s)(g'P) = g \cdot s(g^{-1}g'P) =
g \cdot (g^{-1}g', f(g^{-1}g')) =   (g'P, f(g^{-1}g')).$$ The
vector space $V_\lambda$ of global sections is an irreducible
representation of $G$; the action of $g \in G$ on $s_f$ is $(g
\cdot s_f)(g'P) = g \cdot s_f(g^{-1}g'P)$.

The $N$-th tensor power $L_\lambda^{\otimes N}$ is isomorphic to
$L_{N \lambda}$.

\subsubsection{Choosing a $T$-linearization of $L_\lambda$}

There is a canonical $T$-linearization of $L_\lambda$, given by
restricting the action of $G$ on $L_\lambda$ to $T$.  We shall
call this the ``democratic'' linearization. A weight $\mu$ may be
used to twist the democratic linearization;
$$t \cdot (g,z) = (tg, \mu(t)z).$$ We shall call this the
``$\mu$--linearization''. Indeed the set of all $T$-linearizations
are given by the characters $\mu$ of $T$.

The $\mu$-twisted action of $T$ on a section $s_f$ is given by the
formula,
$$(t \cdot s_f)(gP) = (gP, e^\mu(t)f(t^{-1}g).$$ Hence $s_f$ is
$T$-invariant iff $\mu(t)f(t^{-1}g) = f(g)$ for all $t \in T$.
Equivalently, we have that $s_f$ is $T$-invariant iff for all $t
\in T$ and $g \in G$,
$$f(tg) = e^\mu(t)f(g).$$  The action on a section $s_f$ of
$L_{\lambda}^{\otimes N}$ is given by $(t \cdot s_f)(gP) = (gP,
e^{N \mu}(t)f(g))$, and so the $T$-invariant sections $s_f$ of the
$N$-th tensor power of $L_\lambda$ are those which satisfy
$$f(t \cdot g) = e^{N \mu}(t) f(g).$$

\subsection{The G.I.T. construction}

The G.I.T. quotient $F \q T$ associated to the pair
$(\lambda,\mu)$ is the projective variety,
$$F \q T = \mathrm{Proj}\Bigg(\bigoplus_{N=0}^\infty \Gamma(F, L_{\lambda}^{\otimes
N})^T\Bigg),$$ where $T$ acts on $L_\lambda$ via the
$\mu$-linearization.

\begin{definition}
The set of semistable points $F^{ss} \subset F$ is defined by $p
\in F^{ss}$ iff there exists some positive integer $N$ and a
$T$-invariant global section $s$ of $L_\lambda^{\otimes N}$ such
that $s(p) \neq 0$.  (Normally there is the additional requirement
that $X_s = \{p \in F \mid s(p) \neq 0\}$ is affine but this is
automatic since $F$ is a projective variety.)  If we take the $\mu$-linearization 
of $L_\lambda$, then we shall say that a semistable point is $\mu$-semistable.
\end{definition}

A standard result of Mumford's Geometric Invariant Theory is that
the G.I.T. quotient is topologically given as a quotient space of
the open subset of semistable points.  In particular there is a
surjective continuous map $\pi : F^{ss} \to F \q T$, where $\pi(x)
= \pi(y)$ iff the closures of the $T$-orbits of $x$ and $y$ (Zariski
closure) have non-empty intersection in $F^{ss}$; $\overline{T
\cdot x} \cap \overline{T \cdot y} \cap F^{ss} \neq \emptyset$.
The space $F \q T$ has the quotient topology relative to the
surjective map $\pi$.

The proof of the following theorem will be given in section 
\S\ref{saturation}. This theorem allows us to explicitly construct an ample line bundle of
$F \q T$, and to cover $F \q T$ by explicit affine varieties.

\begin{theorem}\label{maintheorem}(Main Theorem)
Suppose that $\lambda$ is a dominant weight and $\mu$ is any weight, 
such that $\lambda - \mu$ lies in the root lattice of $\SL(n,\C)$.  Then 
if $p \in F$ is $\mu$--semistable there is a global $T$-invariant
section $s$ of $L_\lambda$ such that $s(p) \neq 0$.
\end{theorem}

\begin{remark}
If $\lambda - \mu$ is not in the root lattice, then $\Gamma(F,L_\lambda)^T$ is zero.  
In fact, $\Gamma(F,L_\lambda)^T$ is nonzero if and only if $\lambda - \mu$ is in the 
root lattice and $\mu$ lies within the convex hull of the Weyl orbit of $\lambda$.
If $\mu$ does not lie in the convex hull of the Weyl orbit of $\lambda$, 
then $\Gamma(F,L_\lambda^{\otimes N})^T$ is zero for all $N > 0$; in this case 
there are no semistable points in $F$, and the 
quotient $F \q T$ is empty.
\end{remark}

The following is taken from \cite{Dolgachev}, chapter 8.   Let
$s_1,\ldots,s_m$ be a basis of the $T$-invariant sections of
$L_\lambda$ for the $\mu$-linearization.   By theorem
\ref{maintheorem}, the semistable points $F^{ss}$ are covered by
the affine open subsets $X_{s_i}$, where $X_{s_i} = \{x \in F \mid
s_i(x) \neq 0\}$.  Let $Y_i$ be the affine quotient $X_{s_i} \q
T$; the affine coordinate ring of $Y_i$ is
$\mathcal{O}(X_{s_i})^T$.  The $Y_i$'s may be glued together via
the transition functions $s_i/s_j$ to form the G.I.T. quotient $F
\q T$, and simultaneously an ample line bundle $M$ of $F \q T$,
such that $\pi^\ast(M) = \iota^\ast(L_\lambda)$, where $\iota :
F^{ss} \to F$ is the inclusion map.

As stated in the introduction, it remains an open problem to
compute the minimal integer $N$ such that $M^{\otimes N}$ is very
ample.

\section{Matroid polytopes and weight polytopes}

A matroid is a pair $M = (E, \mathcal{B})$ where $E$ is a finite
set called the ground set of $M$, and $\mathcal{B}$ is a nonempty
collection of subsets of $E$ called bases that satisfy the
exchange condition, which is that for any two bases $B_1,B_2 \in
\mathcal{B}$, if $x \in B_1 \setminus B_2$ then there is an
element $y \in B_2 \setminus B_1$ such that $(B_1 \setminus \{x\})
\cup \{y\} \in \mathcal{B}$ is a basis.  Necessarily it follows
that all bases $B \in \mathcal{B}$ have the same cardinality,
which is called the rank of $M$.  Matroids are a generalization of
finite configurations of vectors, where the only data known about
the set of vectors is which subsets are maximal independent
subsets.  The collection of maximal independent subsets satisfies
the exchange condition.  Similarly, a linear subspace $\Lambda$ of
dimension $k$ of $\C^n$ determines a matroid $M(\Lambda)$, given
by the vector configuration
$\{\pi_\Lambda(e_1),\ldots,\pi_\Lambda(e_n)\}$ where $\pi_\Lambda$
is orthogonal projection onto $\Lambda$ (for the standard
Hermitian form) and the $e_i$'s are the standard basis vectors of
$\C^n$.

\subsection{Matroid polytopes}

Suppose that $M = (E, \mathcal{B})$ is a matroid, and $E =
\{1,2,3,\ldots,n\}$. For each $B \in \mathcal{B}$ let $v^B \in \R^n
/ \Delta$ ($\Delta$ is the diagonal in $\R^n$) be given by $v^B_i =
0$ if $i \notin B$ and $v^B_i = 1$ if $i \in B$.  Let $P_M$ be the
convex hull of $\{v^B \mid B \in \mathcal{B}\}$.  We call $P_M$ a
matroid polytope. Each $v^B$ is a vertex of $P_M$ and so $M$ may be
recovered from $P_M$.

\begin{theorem} (Gel'fand Goresky MacPherson Serganova
\cite{GGMS}) Two vertices $v^{B_1}$, $v^{B_2}$ of $P_M$ form an
edge of $P_M$ iff $v^{B_1} - v^{B_2} = e_i - e_j$ for some $i \neq
j$, where $e_1,\ldots,e_n$ are the standard basis vectors of
$\R^n$.  In other words, edges of $P_M$ are parallel to roots of
$\SL(n,\C)$.  (In fact, the bases $B_1$ and $B_2$ differ by a
single exchange iff $v^{B_1}$ and $v^{B_2}$ form an edge of
$P_M$.)

Conversely, if $P$ is a polytope where all vertices are $0/1$
vectors (each component is either $0$ or $1$), and each edge of
$P$ is parallel to an $\SL(n,\C)$ root, then there is a matroid $M$ such that
$P = P_M$.
\end{theorem}

\begin{remark}
A natural way that matroid polytopes arise is by restricting the
momentum mapping $\rho : \Gr_k(\C^n) \to \mathfrak{t}^\ast_0$ for
the action of the maximal compact subtorus $T_0$ in $T$ on the
Grassmannian to the closure of an orbit $T \cdot \Lambda$, see
\cite{GGMS} or \cite{BorovikGelfandWhite}. The polytope
$P_{M(\Lambda)}$ is the image of $\rho$ restricted to the closure of
$T \cdot \Lambda$.
\end{remark}

\subsection{Weight polytopes}

Suppose that $V$ is a finite dimensional complex representation of
a torus $T$. Then $V$ is a direct sum of weight spaces,
$$V = \bigoplus_{\mu} V[\mu],$$
where $V[\mu] = \{v \in V \mid \text{$t \cdot v = e^\mu(t)v$ for all
$t \in T$}\}$.  Note that a section $s \in V_\lambda = \Gamma(F,
L_\lambda)$ is $T$-invariant for the $\mu$-linearization of
$L_\lambda$ if and only if $s \in V_\lambda[\mu]$.

Given a dominant weight $\lambda$ let $P_\lambda$ denote the
associated parabolic subgroup.  For each $g \in G$, let
$$wt_\lambda(g) = \{\mu \mid (\exists s \in
V_\lambda[\mu])(s(gP_\lambda) \neq 0)\}.$$ Let the \emph{weight
polytope} $\overline{wt}_\lambda(g)$ be the convex hull of
$wt_\lambda(g)$ (the convex hull is taken in $\mathfrak{t}_0^\ast$,
where $\mathfrak{t}_0$ is the Lie algebra of the maximal compact
torus $T_0 \subset T$).

\begin{lemma}
For any two dominant weights $\lambda_1$ and $\lambda_2$, we have
$$wt_{\lambda_1}(g) + wt_{\lambda_2}(g) = wt_{\lambda_1 +
\lambda_2}(g),$$ where the summation denotes the Minkowski sum, $A
+ B = \{a+b \mid a \in A, b \in B\}$.
\end{lemma}

\begin{proof}
Suppose that $\mu_1 \in wt_{\lambda_1}(g)$ and $\mu_2 \in
wt_{\lambda_2}(g)$.  Let $s_1 \in V_{\lambda_1}[\mu_1]$ and $s_2
\in V_{\lambda_2}[\mu_2]$ such that $s_1(gP_{\lambda_1}) \neq 0$
and $s_2(gP_{\lambda_2}) \neq 0$.  Recall there are functions
$f_1, f_2:G \to \C$ such that $s_1 = s_{f_1}$ and $s_2 = s_{f_2}$.
We have that $f_1(g) \neq 0$ and $f_2(g) \neq 0$.  Hence, $f_1(g)
f_2(g) \neq 0$.  The section $s_{f_1 f_2}$ lies in $V_{\lambda_1 +
\lambda_2}[\mu_1 + \mu_2]$, and is nonzero at $gP_{\lambda_1 +
\lambda_2}$.

Now suppose that $\mu \in wt_{\lambda_1 + \lambda_2}(g)$. We may
identify the irreducible representation $V_\lambda$ as the space
of global sections of $\pi^\ast(L_\lambda)$ of $G/B$ where $B$ is
the Borel subgroup of $G$ and $\pi : G/B \to G/P_{\lambda}$.  This
is justified since the pullback $\pi^\ast : \Gamma(G/P_{\lambda},
L_{\lambda}) \to \Gamma(G/B, \pi^\ast(L_{\lambda}))$ is an
isomorphism of vector spaces.  We shall also abuse notation and
identify $L_\lambda$ with the pullback $\pi^\ast{L_\lambda}$.

The tensor product $V_{\lambda_1} \otimes V_{\lambda_2}$ is the
vector space of sections of the outer tensor product
$L_{\lambda_1} \boxtimes L_{\lambda_2}$ of $G/B \times G/B$, where
$B$ is the Borel subgroup.   The irreducible representation
$V_{\lambda_1 + \lambda_2}$ is a direct summand of $V_{\lambda_1}
\otimes V_{\lambda_2}$, and the projection $V_{\lambda_1} \otimes
V_{\lambda_2} \to V_{\lambda_1 + \lambda_2}$ is realized by
pulling back $L_{\lambda_1} \boxtimes L_{\lambda_2}$ to the
diagonal $\Delta \subset G/B \times G/B$.  We have assumed there
is a section $s \in V_{\lambda_1 + \lambda_2}[\mu]$ such that
$s(gB) \neq 0$. Clearly $(V_{\lambda_1} \otimes
V_{\lambda_2})[\mu]$ surjects onto $V_{\lambda_1 +
\lambda_2}[\mu]$.  Furthermore,
$$(V_{\lambda_1} \otimes V_{\lambda_2})[\mu] = \sum_{\mu_1 + \mu_2 =
\mu} V_{\lambda_1}[\mu_1] \otimes V_{\lambda_2}[\mu_2].$$ Hence
there must exist weights $\mu_1, \mu_2$ such that $\mu_1 + \mu_2 =
\mu$ and some component $s' = s_1 s_2$ of $s$ such that
$s_1(gB)s_2(gB) \neq 0$ and $s_1 \in V_{\lambda_1}[\mu_1]$ and
$s_2 \in V_{\lambda_2}[\mu_2]$.
\end{proof}

\begin{corollary}
Suppose that $\lambda = \sum_{i=1}^{n-1} a_i \varpi_i$ is dominant,
i.e. each $a_i$ is non-negative and $\varpi_i$ denotes the $i$-th
fundamental weight connected to the Grassmannian $\Gr_i(\C^n)$. Then
for any $g \in G$,
$$wt_\lambda(g) = \sum_{i=1}^{n-1} a_i \cdot wt_{\varpi_i}(g),$$
where the sum indicates Minkowski sum and $a_i \cdot
wt_{\varpi_i}(g)$ denotes the $a_i$-fold Minkowski sum of
$wt_{\varpi_i}(g)$.
\end{corollary}

The weight polytope $\overline{wt}_\lambda(g)$ is a \emph{flag
matroid polytope} within the more general setting of Coxeter matroid
polytopes, see \cite{BorovikGelfandWhite}.  However, we will only
need to consider standard matroid polytopes, as they are the
building blocks for flag matroid polytopes.

\begin{proposition}
Suppose that $\varpi_k$ is the $k$-th fundamental weight. Then
$\overline{wt}_{\varpi_k}(g)$ is a matroid polytope for any $g \in
G$.
\end{proposition}

\begin{proof}
A basis for the sections of $L_{\varpi_k}$ is given by bracket
functions $[i_1,i_2,...,i_k]$ where $1 \leq i_1 < i_2 < \cdots <
i_k \leq n$.  The section $s=[i_1,i_2,...,i_k]$ is equal to $s_f$,
where $f : G \to \C$ assigns the determinant of the $k$ by $k$
submatrix given by columns $1,2,\ldots,k$ and rows
$i_1,i_2,...,i_k$ of $g \in G$.  The bracket $[i_1,i_2,...,i_k]$
belongs to the weight space $V_{\varpi_k}[\mu]$ where $e^\mu =
(a_1,a_2,\ldots,a_n) \in \Z^n / \Delta$ is given by $a_i = 1$ if
$i = i_t$ for some $t$, $1 \leq t \leq k$, otherwise $a_i = 0$.
Now suppose that $gP_{\varpi_k} \in G/P_{\varpi_k} = \Gr_k(\C^n)$.
The linear subspace defined by $gP_{\varpi_k}$ is the span of the
first $k$ columns of $g$.  We have that $\mu \in
\overline{wt}_{\varpi_k}(g)$ iff $\mu$ is a $0/1$ vector (mod
$\Delta$) with $k$ ones (occuring at $I = (i_1,i_2,\ldots,i_k)$)
and $n-k$ zeros such that the $I$-th minor is nonzero.

Let $M(g)$ be the matroid with ground set $\{1,2,\ldots,n\}$ of
the vector configuration $r_1,r_2,\ldots,r_n \in \C^k$ where $r_i$
is the $i$-th row of $g$ restricted to the first $k$ columns, i.e.
$r_i = (g_{i,1},g_{i,2},\ldots,g_{i,k})$.  It is clear that the
matroid polytope of $M(g)$ is the weight polytope
$\overline{wt}_{\varpi_k}(g)$.
\end{proof}

\section{Saturation properties of weight polytopes}\label{saturation}

We shall prove the following lemma by a combinatorial argument.
The main theorem \ref{maintheorem} easily follows from this lemma.
Neil White proved in \cite{White} the exact same statement for
$\lambda = \varpi_k$, using a theorem of Edmonds in matroid
theory.

\begin{lemma}\label{saturationlemma}  Suppose $g \in G$ and $\lambda$ is a dominant
weight. Suppose $\mu$ is a weight such that $\lambda - \mu$ is in the root lattice. 
Then $N \mu \in wt_{N
\lambda}(g)$ implies $\mu \in wt_\lambda(g)$ for all $N > 0$.
\end{lemma}

\begin{remark}
If $G$ is any complex semi-simple group, and $\lambda$ is a dominant
weight, and $\lambda - \mu$ is in the root lattice of $G$, 
then $V_{N \lambda}[N \mu] \neq 0$ implies $V_{\lambda}[\mu]
\neq 0$. However, the lemma is a much stronger statement than this
(and it does not hold for groups other than $\SL(n,\C)$) because $g$
is fixed (i.e., the point $gP_\lambda \in G/P_\lambda$ is fixed).
\end{remark}

Let $R$ be the set of $\SL(n,\C)$ roots.  Let $Q(R)$ (resp. $P(R)$)
denote the root lattice (resp. weight lattice).  Convex hulls of
subsets of the weight lattice, denoted by an overline, should take
place in $\mathfrak{t}_0^\ast$, which is isomorphic to $P(R) \otimes
\R = \overline{P(R)} = \R^n / \Delta$. The map $\epsilon : P(R) \to
\Z / n\Z$ given by $\epsilon(a_1,\ldots,a_n) = \sum_i a_i \mod n$ is
a homomorphism of abelian groups, and $Q(R) = \ker(\epsilon)$.

\begin{definition}
A finite subset $A$ of $Q(R)$ is called {\it root--saturated} if
\begin{itemize}
\item the convex hull $\overline{A}$ is such that each edge $e_i$ is
parallel to a root $\gamma_i$ in $R$, (i.e. $\overline{A}$ is a flag
matroid polytope, see \cite{BorovikGelfandWhite}.)
\item for each $x \in A$, $A = \overline{A} \cap Q(R)$.
\end{itemize}
\end{definition}

We will eventually prove that $wt_\lambda(g) -\lambda$ (the set
$wt_\lambda(g)$ translated by $-\lambda$) is root-saturated for any
dominant weight $\lambda$.

\begin{definition}
The {\it Minkowski sum} of two subsets $A$,$B$ of Euclidean space, denoted $A+B$, is
$$A+B = \{a+b \,:\, a \in A,\, b \in B\}$$
\end{definition}

\begin{lemma}
Suppose that $\alpha_1,\ldots,\alpha_{n-1} \in R$ are independent
over $\Q$.  Then they are a basis for the root lattice $Q(R)$.
\end{lemma}

\begin{proof}
The proof goes by induction on $n$.  If $n=2$ there are only two roots $\alpha, -\alpha$ and they generate the
same lattice.  Now suppose that $n > 2$.
Let $\Z[\alpha_1,\ldots,\alpha_{n-1}]$ be the $\Z$--span of $\alpha_1,\ldots,\alpha_{n-1}$.
Without loss of generality we may assume that each $\alpha_i$ is a positive root since negating $\alpha_i$ does not change the span over $\Z$.  Let $\sigma_1,\ldots,\sigma_{n-1}$ be the standard simple roots of $SL(n)$.  That is, $\sigma_i = e_i - e_{i+1}$.  Note that any positive root $e_i - e_j = \sum_{t=i}^{j-1} \sigma_t$ is a sum of consecutive simple roots.  Conversely any consecutive sum of simple roots is a positive root.
We may choose some $w \in W$ (where $W$ is the Weyl group) such that $w(\alpha_{n-1}) = \sigma_{n-1}$.
In particular if $\alpha_{n-1} = e_i - e_j$ let $w$ be the product of two cycles $(n-1 \quad i)(n \quad j)$.
Since elements of $W$ induce isomorphisms of the lattice $Q(R)$, we have that $w(\alpha_1),\ldots,w(\alpha_{n-1})$ is a basis of $Q(R)$ if and only if $\alpha_1,\ldots,\alpha_{n-1}$ is a basis of $Q(R)$.
Reassign $\alpha_i := w(\alpha_i)$.
For each $i \leq n-2$, if $\alpha_i = e_s - e_n = \sigma_s + \cdots + \sigma_{n-1}$ replace
$\alpha_i$ with $\alpha_i - \sigma_{n-1} = \sigma_s + \cdots + \sigma_{n-2} = e_s - e_{n-1}$.
Now the roots $\alpha_1,\ldots,\alpha_{n-2}$ may be identified with roots of $SL(n-1)$.
By the induction hypothesis $\Z[\alpha_1,\ldots,\alpha_{n-2}] =  \Z[\sigma_1,\ldots,\sigma_{n-2}]$.  Since $\alpha_{n-1} = \sigma_{n-1}$ we have that
$\Z[\alpha_1,\ldots,\alpha_{n-1}] = Q(R)$.

\end{proof}

\begin{lemma}
Suppose that $A$ and $B$ are root-saturated, and $\overline{A} \cap \overline{B}$ is nonempty.
Then $A \cap B$ is nonempty.
\end{lemma}

\begin{proof}
The proof is by induction on the dimension of $\overline{A}$.  If
$\dim \overline{A} = 0$ then $A = \{a\}$ for some $a \in Q(R)$.
Then $\overline{A} \cap \overline{B} = A \cap B = \{a\}$. Now
suppose that $\dim \overline{A} \geq 1$.

We have two cases, the first case is that the intersection
$\overline{A} \cap \overline{B}$ contains a boundary point of
$\overline{A}$. Then there is some facet $F$ of $\overline{A}$ such
that $F \cap \overline{B}$ is nonempty. We claim $F \cap A$ is
root--saturated.  The vertices of $F$ are within $F \cap A$. so
$\overline{F \cap A} \supset F$.  On the other hand $F \subset F
\cap A$ so $F \subset \overline{F \cap A}$; therefore $F =
\overline{F \cap A}$. The edges of $F$ are also edges of
$\overline{A}$ hence they are parallel to roots. Furthermore, for
any $x \in F \cap A$, we have $\overline{F \cap A} \cap Q(R) \subset
A$ since $A$ is root--saturated, and it follows that $\overline{F
\cap A} \cap Q(R) = F \cap A$ since $F \cap A \subset A \subset
Q(R)$. Since $\dim F < \dim \overline{A}$ we may apply the induction
hypothesis to get that $F \cap A \cap B$ is nonempty and hence $A
\cap B$ is nonempty.

On the other hand suppose that $\overline{A} \cap \overline{B}$
contains no boundary point of $\overline{A}$.  Let $L_A(R)$ be the sub--lattice of $Q(R)$ spanned by
the roots which are parallel to some edge of $\overline{A}$.
Let $a_0 \in A$ be a vertex of $\overline{A}$.  Note that the affine space
$H_A = a_0 + \overline{L_A(R)}$ is the smallest affine space containing $\overline{A}$.  We claim
$H_A \cap \overline{B} = \overline{A} \cap \overline{B}$.
Suppose that $z \in H_A \cap \overline{B}$.  Let $a \in \overline{A} \cap \overline{B}$.
Since $H_A$ has the same dimension as $\overline{A}$, there are linear inequalities
$\eta_i(x) \leq f_i$ where the interior of $\overline{A}$ consists of points $x \in H_A$
where the inequalities are strict; that is, $\eta_i(x) < f_i$ for all $i$ if and only if $x$ is an
interior point of $\overline{A}$.  The boundary points of $\overline{A}$ are those
points $x \in \overline{A}$ such that $\eta_i(x) = f_i$ for some $i$.
Let $c(t) = (1-t)a + t z$ for $0 \leq t \leq 1$.  Suppose that $z \notin \overline{A}$.
Then there is some $i$ such that $\eta_i(z) > f_i$.  However $a$ is an interior point of $\overline{A}$ and
so $\eta_i(a) < f_i$.  Hence there is some $t_0$ such that $\eta(c(t_0)) = f_i$ in which case
$c(t_0)$ is a boundary point of $\overline{A}$.  But $c(t) \in \overline{B}$ for each $t$ by convexity
of $\overline{B}$.  This contradicts that $\overline{A} \cap \overline{B}$
is disjoint from the boundary of $A$.  Hence $H_A \cap \overline{B} = \overline{A} \cap \overline{B}$.
Therefore $(H_A \cap Q(R)) \cap B = A \cap B$ since $\overline{A} \cap Q(R) = A$ and $\overline{B} \cap Q(R) = B$.

We now show by induction on $\dim \overline{B}$, that for any $B$ which is root-saturated, that
$H_A \cap \overline{B}$ is nonempty implies $(H_A \cap Q(R)) \cap B$ is nonempty.
Suppose that $\dim \overline{B} = 0$.  Then $B = \{b\}$ for some $b \in Q(R)$, and so
$b \in (H_A \cap Q(R)) \cap B$.  Now suppose that $\dim \overline{B} \geq 1$.  We have two cases.

First suppose that $H_A \cap \overline{B}$ intersects the boundary of $\overline{B}$ nontrivially.
Then there is a face $F$ of $\overline{B}$ such that $H_A \cap F$ is nonempty.  Since
$F \cap B$ is root--saturated, $\overline{F \cap B} = F$, $H_A \cap F$ is nonempty, and $\dim F < \dim B$, we may apply the induction hypothesis and we're finished.

Now suppose that $H_A \cap \overline{B}$ is disjoint from the boundary of $\overline{B}$.
Let $L_B(R)$ be the sub--lattice of $Q(R)$ spanned by the roots which are parallel to some edge of
$\overline{B}$.  Let $b_0 \in B$ be a vertex of $\overline{B}$.
The affine space $H_B = b_0 + \overline{L_B(R)}$ is the smallest affine space containing $\overline{B}$.
As above, we have that $H_A \cap H_B = H_A \cap \overline{B}$ and so
$(H_A \cap Q(R)) \cap (H_B \cap Q(R)) = A \cap B$.  Since $H_A$ does not intersect
the boundary of $\overline{B}$, we have that $H_A \cap H_B$ is a single point $z_0$, since if the dimension of
the intersection $H_A \cap H_B = H_A \cap \overline{B}$ is greater than zero then $H_A \cap \overline{B}$ is unbounded.  But $\overline{B}$ is compact since $B$ is finite and this cannot happen.
We now show that $z_0 \in Q(R)$.  We have that $z_0 = a_0 + v_A = b_0 + v_B$ where $a_0 \in A$, $b_0 \in B$,
$v_A \in \overline{L_A(R)}$, $v_B \in \overline{L_B(R)}$.  Let $\{\alpha_1,\ldots,\alpha_p\}\subset R$
be a basis of $L_A(R)$ and let $\{\beta_1,\ldots,\beta_q\} \subset R$ be a basis of $L_B(R)$.
Since the intersection of $H_A$ and $H_B$ is a point, we have that
$\overline{L_A(R)} \cap \overline{L_B(R)} = \{0\}$.
Hence the set $\{\alpha_1,\ldots,\alpha_p,\beta_1,\ldots,\beta_q\}$ is linearly independent in
$\overline{Q(R)}$.  Choose $\{\gamma_1,\ldots,\gamma_r\} \subset R$ so that
$\{\alpha_1,\ldots,\alpha_p,\beta_1,\ldots,\beta_q,\gamma_1,\ldots,\gamma_r\}$ is a basis for $\overline{Q(R)}$.
By the Lemma above this is also a basis for the lattice $Q(R)$.
Now $v_A = \sum_i c_i \alpha_i$ and $v_B = \sum_j d_j \beta_j$ are unique expressions for $v_A,v_B$.
But also the difference $a_0 - b_0 = v_B - v_A = (\sum_j d_j \beta_j) - (\sum_i c_i \alpha_i)$ lies within
the lattice $Q(R)$, and so the coefficients $c_i$,$d_j$ must be integers.
Hence, $z_0$ is a lattice point and we've finished the proof of the Lemma.
\end{proof}

\begin{theorem}\label{SaturationSumTheorem}
Suppose that $A$ and $B$ are root-saturated.  Then the Minkowski
sum $A + B$ is root-saturated.
\end{theorem}

\begin{proof}
We show that the Minkowski sum $A+B$ is root-saturated if $A$ and
$B$ are each root-saturated. Clearly $A+B$ is finite, and the
elements are within $Q(R)$ since $Q(R)$ is closed under addition.
First we show that the edges of $\overline{A+B}$ are parallel to
roots.  Clearly $\overline{A} + \overline{B} = \overline{A+B}$.
The Minkowski sum of two polytopes $P,Q$ has edges of the
following types:
\begin{itemize}
\item (vertex of P) + (edge of Q).
\item (edge of P) + (vertex of Q).
\item (edge of P) + (edge of Q), providing these edges are parallel.
\end{itemize}
We leave the proof to the reader (the proof is easily obtained by
observing that the fan of $P+Q$ is the meet of the fan of $P$ with
the fan of $Q$).  In all three cases, the resulting edge is parallel
to an edge of either $P$ or $Q$ or both, and hence it is parallel to
some root in $R$.

Next we must show that $A+B = (\overline{A+B}) \cap Q(R)$.  Suppose that $z \in (\overline{A+B}) \cap Q(R)$.
Hence there exists $x \in \overline{A}$ and $y \in \overline{B}$ such that $x+y = z$.  Hence
$x \in (z + \overline{-B}) \cap \overline{A}$, where $-B = \{-b : b \in B\}$.  Clearly $z + (-B)$ is root-saturated.
Hence, we may apply the Lemma above to get a lattice point $x_0$ in the intersection.  Since $A$ is
saturated, we have that $x_0 \in A$.  Now we have that $z = x_0 + y_0$ where $y_0 \in \overline{B}$.
But since $z,x_0 \in Q(R)$ we have that $y_0 = z - x_0 \in Q(R)$, and so $y_0 \in B$ since $B$ is root-saturated, and we're finished.
\end{proof}

\begin{lemma}
If $\varpi_k$ is a fundamental weight and $g \in G$ then the
translation $wt_{\varpi_k}(g) - \varpi_k$ is root-saturated.
\end{lemma}

\begin{proof}
Note that all elements of $wt_{\varpi_k}(g)$ are $0/1$ vectors (mod
$\Delta$) having $k$ ones and $n-k$ zeros.  Translating by $-
\varpi_k$ results in vectors whose first $k$ components may be
either $0$ or $-1$ and last $n-k$ components are $0$ or $+1$, and
the sum of all components is zero.  Hence the first $k$ components
define a vertex of the negated unit $k$-cube, and the last $n-k$
components are vertices of the $n-k$-cube.  Therefore, there can be
no additional lattice points in the convex hull. We already showed
that the convex hull of $wt_{\varpi_k}(g)$ is a matroid polytope, so
the edges are parallel to roots.  This property is preserved by
translations.
\end{proof}

\begin{corollary}\label{dom_saturated}
For any dominant weight $\lambda$ and $g \in G$, the set
$wt_\lambda(g) - \lambda$ is root-saturated.
\end{corollary}

\begin{proof}
We have that $\lambda = \sum_{k=1}^{n-1} a_k \varpi_k$, where the
$a_k$'s are non-negative integers.  Also, $wt_\lambda(g) =
\sum_{k=1}^{n-1} a_k \cdot wt_{\varpi_k}(g)$ (Minkowski sum). Hence,
$$wt_\lambda(g)- \lambda = \sum_{k=1}^{n-1} a_k \cdot
(wt_{\varpi_k}(g) - \varpi_k).$$ Since the root-saturated property
is preserved under Minkowski sums, we have that $wt_\lambda(g) -
\lambda$ is root-saturated.
\end{proof}

{\it Proof of lemma \ref{saturationlemma}.} \begin{proof} Suppose
that $N \mu \in wt_{N \lambda}(g)$.  Then $N(\mu - \lambda) \in
wt_{N\lambda}(g) - N\lambda$.  The convex hull of $wt_{N\lambda}(g)
- N\lambda$ scaled by $1/N$ is equal to the convex hull of
$wt_{\lambda}(g) - \lambda$ since $N \cdot wt_\lambda(g) = wt_{N
\lambda}(g)$. Therefore $\mu - \lambda$ is in the convex hull of
$wt_\lambda(g) - \lambda$. But since $\mu - \lambda \in Q(R)$ and
$wt_\lambda(g) - \lambda$ is root-saturated, we have that $\mu -
\lambda \in wt_\lambda(g) - \lambda$, so $\mu \in wt_\lambda(g)$.
\end{proof}

{\it Proof of main theorem \ref{maintheorem}.} \begin{proof}
Suppose that $gP_\lambda$ is semistable relative to the
$\mu$-linearization of the line bundle $L_\lambda$.  This means
there is some $N > 0$ and a section $s \in \Gamma(G/P_\lambda,
L_\lambda^{\otimes N})^T$ such that $s(gP_\lambda) \neq 0$.  This
means that $N \mu \in wt_{N \lambda}(g)$.  By Lemma
\ref{saturationlemma} we have that $\mu \in wt_\lambda(g)$.  So
there must exist a section $s' \in \Gamma(G/P_\lambda,
L_\lambda)^T$ such that $s'(gP_\lambda) \neq 0$.
\end{proof}

\subsection{Failure of main theorem for $G = \mathrm{SO}(5,\C)$}

Let $B(z,w)$ be the bilinear form on $\C^5$ given by
$$B(z,w) = z_1 w_5 + z_2 w_4 + z_3 w_3 + z_4 w_2 + z_5 w_1 = 2 z_1
w_5 + 2 z_2 w_4 + z_3 w_3.$$ Now $\mathrm{SO}(5,\C) \subset
\SL(5,\C)$ is the subgroup preserving $B$.  The maximal torus $T$
may be taken to the diagonal matrices in $\mathrm{SO}(5,\C)$.
Elements of $T$ have the form $diag(t_1,t_2,1,1/t_2,1/t_1)$ for
$t_1,t_2 \in \C^\ast$. Let $\varpi_1$ denote the first fundamental
weight of $\mathrm{SO}(5,\C)$.  We have
$e^{\varpi_1}(t_1,t_2,1,1/t_1,1/t_2) = t_1$, but the second
fundamental weight does not lift to a character of
$\mathrm{SO}(5,\C)$ - one needs to go the universal cover to find
such a character. Let $P_{\varpi_1} \subset \mathrm{SO}(5,\C)$ be
the associated parabolic subgroup. The quotient space
$\mathrm{SO}(5,\C)/P_{\varpi_1}$ may be identified with the space
of isotropic lines in $\C^5$.

Let $x$ be the (isotropic) line through
$(1,\sqrt{-1},0,\sqrt{-1},1)$.  Let $g_x \in \mathrm{SO}(5,\C)$ be
such that $g_x P_{\varpi_1} = x$. The set $wt_{\varpi_1}(g_x)$ is
equal to $\{\varpi_1, 2\varpi_2-\varpi_1, -2\varpi_2 + \varpi_1,
-\varpi_1\}$. Depiction:
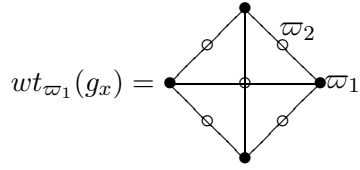
\[wt_{\varpi_1}(g_x) =
\begin{xy}
(0,-10)*{}="A"; (10,0)*{}="B"; (0,10)*{}="C"; (-10,0)*{}="D";
(0,0)*{}="O"; "A"; "B"; **\dir{-}; "A"; "O"; **\dir{-}; "B"; "C";
**\dir{-}; "C"; "D"; **\dir{-}; "A"; "D"; **\dir{-}; "B"; "O";
**\dir{-}; "C"; "O"; **\dir{-}; "D"; "O"; **\dir{-};
"A"*{\bullet}; "B"*{\bullet}; "C"*{\bullet}; "D"*{\bullet};
"O"*{\circ}; (5,5)*{\circ}; (5,-5)*{\circ}; (-5,5)*{\circ};
(-5,-5)*{\circ}; (7,7)*{\varpi_2}; (13,0)*{\varpi_1}
\end{xy}
\]
This set is missing the origin, although $V_{\varpi_1}[0] \neq 0$
and $\varpi_1 \in Q(\mathrm{SO}(5,\C))$, so $wt_{\varpi_1}(g_x) -
\varpi_1$ is not root-saturated. Note the origin does belong to
$wt_{2 \varpi_1}(g_x) = wt_{\varpi_1}(g_x) + wt_{\varpi_1}(g_x)$.
Therefore $x$ is semistable for the democratic linearization of
$L_{\varpi_1}$.  It follows that for the democratic linearization
of $L_{\varpi_1}$, one requires a $T$-invariant section of
$L_{\varpi_1}^{\otimes 2}$ to pick out the semistable point $x$.

\section{Projective Normality}

Let $H$ be the group of diagonal matrices in $\GL(n,\C)$.  Hence
$T \subset H$ is the set of diagonal matrices with determinant
one. Let $\chi_1,\ldots,\chi_m$ be $m$ characters of $H$.  That
is, each $\chi_i : H \to \C^\ast$ is an algebraic homomorphism of
groups. Each $\chi_i$ is given by a point $\ba_i =
(a_{i,1},\ldots,a_{i,n}) \in \Z^n$, where
$$\chi_i(h_1,\ldots,h_n) = \prod_{j = 1}^n h_j^{a_{i,j}}.$$  These
characters determine an action of $H$ on $\mathbb{A}^m$ by
$$h \cdot (z_1,z_2,\ldots,z_m) = (\chi_1(h) z_1, \chi_2(h) z_2,
\ldots, \chi_m(h) z_m).$$ Now take any point $z \in \mathbb{A}^m$,
and let $X(z)$ be the Zariski closure of the $H$-orbit of $z$.  That
is, $X(z) = cl(H \cdot z)$.  Certainly $X(z)$ contains a dense torus
and there is a natural action of this torus on $X(z)$; so $X(z)$ is
a (possibly non-normal) toric variety.

But when is $X(z)$ a \emph{normal} toric variety, i.e. when is 
the coordinate ring of $X(z)$ integrally closed in its field of fractions?  
Some notation: if
$A$ is a finite subset of $\Z^d$ then let $\Z(A)$ be the sub-lattice
generated by $A$, let $\N(A)$ be the semigroup of all non-negative
integral combinations of elements of $A$, and let $\Q_0^+(A)$ be the
rational cone in $\Q^d$ given by all non-negative rational
combinations of elements of $A$.  According to Proposition 13.5 of 
\cite{Sturmfels} we have that the semigroup algebra $\C[\N(A)]$ is normal iff 
$\N(A) = \Z(A) \cap \Q_0^+(A)$.

The following proposition is likely well known but we give a
proof for lack of reference.

\begin{proposition}\label{orbitprop}
Let $\supp(z) = \{i \mid z_i \neq 0\}$.  Let $A(z) = \{\chi_i \mid i
\in \supp(z)\}$.  Then $X(z)$ is isomorphic to the affine toric
variety defined by $A(z) \subset \Z^n$.  That is, $X(z)$ is
isomorphic to the affine variety $V \subset \C^{\# A(n)}$ of the
semigroup algebra $\C[\N(A(z))]$, where $\N(A(z))$ is the semigroup
in $\Z^n$ generated by $A(z)$.  Hence $X(z)$ is normal if and only
if $\Z(A(z)) \cap \Q_0^+(A(z)) = \N(A(z))$.
\end{proposition}

\begin{proof}
Let $\bar{z} \in \C^m$ be given by $\bar{z}_i = 1$ if $i \in
\supp(z)$ and $\bar{z}_i = 0$ otherwise.  Let $s_i = 1/z_i$ if $z_i
\neq 0$ and $s_i = 1$ if $z_i = 0$.  Then the matrix
$diag(s_1,\ldots,s_m)$ defines an algebraic automorphism of
$\mathbb{A}^m$ which takes $X(z)$ to $X(\bar{z})$, so $X(\bar{z})$
is isomorphic to $X(z)$.  Hence we may assume that all components of
$z$ are either $0$ or $1$.  Additionally, $X(z)$ lives entirely
within the components $i$ where $z_i$ is nonzero.  Hence, we may
project $X(z)$ onto the linear subspace given by the components in
$\supp(z)$, which defines an isomorphism of $X(z)$ onto its image.
Thus, we may assume that each component of $z$ is equal to one. If
$\chi_i = \chi_j$ for some $i,j$, we may also project away one of
these.  Hence we have reduced to the case that the $\chi_i$'s are
distinct, and $z$ is the vector of all ones.   The coordinate ring
of $X(z)$ is now easily seen to be the semigroup algebra
$\C[\N(A(z))]$.
\end{proof}

A dominant weight $\lambda$ of $\SL(n,\C)$ may be lifted to a
dominant weight $\widetilde{\lambda}$ of $\GL(n,\C)$ by
normalizing $\lambda$ so that the last component is zero. That is,
the image of $\widetilde{\lambda} \in \Z^n$ in $\Z^n / \Delta$ is
$\lambda$, and $\widetilde{\lambda}_n = 0$. Let
$$|\widetilde{\lambda}| = \sum_{i=1}^n \widetilde{\lambda}_i.$$
Now $V_\lambda$ is also an irreducible representation of
$\GL(n,\C)$, where $z I_n \in \GL(n,\C)$ acts by scaling each
vector $s \in V_\lambda$ by $z^{|\widetilde{\lambda}|}$, and so if
$\widetilde{g} = z g$ where $z \in \C^\ast$ and $g \in \SL(n,\C)$
then the action of $\widetilde{g}$ is defined by $\widetilde{g}
\cdot s = z^{|\widetilde{\lambda}|} (g \cdot s)$. A basis for the
representation $V_\lambda$ is given by semi-standard tableaux
$\tau$ of shape $\widetilde{\lambda}$ (with total number of slots
equal to $|\widetilde{\lambda}|$), filled with indices from $1$ to
$n$. A section $s_\tau \in V_\lambda[\mu]$ iff the number of times
the index $i$ appears in $\tau$ is equal to $\mu_i$. Here we are
treating $\mu$ as a weight of $\GL(n,\C)$. Note that if
$V_\lambda[\mu] \neq 0$ then $|\mu| = \sum_{i=1}^n \mu_i =
|\widetilde{\lambda}|$ since $|\mu|$ must equal the total number
of slots in $\tau$, where $s_\tau \in V_\lambda[\mu]$.

Recall that $H = \C^\ast(T)$ is the maximal torus in $\GL(n,\C)$
consisting of diagonal matrices.  For each $g \in \GL(n,\C)$ let
$$wt_{\widetilde{\lambda}}(g) = \{\mu  \mid
(\exists s \in V_\lambda[\mu])(s(g P_{\widetilde{\lambda}}) \neq
0)\},$$ where $P_{\widetilde{\lambda}} \subset \GL(n,\C)$ is the
parabolic subgroup $\C^\ast(P_\lambda)$ associated to
$\widetilde{\lambda}$.  Each $\mu \in wt_{\widetilde{\lambda}}(g)
\subset \Z^n$ satisfies $|\mu| = |\widetilde{\lambda}|$.

Note that the root lattice of $\SL(n,\C)$ may be identified with
integral vectors $v \in \Z^n$ whose components sum to zero.
Hence, for any $g \in \SL(n,\C)$ we have an identification of
$wt_\lambda(g) - \lambda$ with $wt_{\widetilde{\lambda}}(g) -
\widetilde{\lambda}$.  In particular, $wt_{\widetilde{\lambda}}(g)
- \widetilde{\lambda}$ is root-saturated.

Let $N_{\widetilde{\lambda}}$ be the sub-lattice of $\Z^n$ given
by
$$N_{\widetilde{\lambda}} = \{v \in \Z^n \mid \; |v| = \sum_{i=1}^n v_i \equiv 0 \!\! \mod |\widetilde{\lambda}|\}.$$

\begin{lemma} For any $g \in \SL(n,\C)$,
$$\Q_0^+(wt_{\widetilde{\lambda}}(g)) \cap N_{\widetilde{\lambda}} = \N(wt_{\widetilde{\lambda}}(g)).$$
\end{lemma}

\begin{proof}
Suppose that $v \in \Q_0^+(wt_\lambda(g)) \cap
N_{\widetilde{\lambda}}$. Then $|v| = d |\widetilde{\lambda}|$ for
some $d \in \N$.  Hence $v$ belongs to the convex hull of the
$d$-th dilate of $wt_{\widetilde{\lambda}}(g)$, so $v$ is in the
convex hull of $wt_{d \widetilde{\lambda}}(g)$, since $wt_{d
\widetilde{\lambda}}(g)$ is the $d$-fold Minkowski sum of
$wt_{\widetilde{\lambda}}(g)$.  But $wt_{d \widetilde{\lambda}}(g)
- d \widetilde{\lambda}$ is root-saturated, and since $v - d
\widetilde{\lambda} \in Q(R)$ we have that $v - d
\widetilde{\lambda} \in wt_{d \widetilde{\lambda}}(g) - d
\widetilde{\lambda}$. Equivalently, $v \in wt_{d
\widetilde{\lambda}}(g)$. Since $wt_{d \widetilde{\lambda}}(g)$ is
the $d$-fold Minkowski sum of $wt_{\widetilde{\lambda}}(g)$, we
have that $v \in \N(wt_{\widetilde{\lambda}}(g))$.
\end{proof}

\begin{corollary}
The semigroup algebra $\C[\N(wt_{\widetilde{\lambda}}(g))]$ is
normal.
\end{corollary}

Now suppose that $\lambda$ is dominant and $P_\lambda$ is the
associated parabolic subgroup.  Choose a basis
$(s_1,s_2,\ldots,s_N)$ of $V_\lambda = \Gamma(\SL(n,\C)/P_\lambda,
L_\lambda)$ such that each basis vector is a generalized
eigenvector for the democratic $T$-action.  (Recall the democratic
action is the restriction of the natural action of $\SL(n,\C)$ on
$V_\lambda$ to $T$.) Let $\iota_\lambda : \SL(n,\C)/P_\lambda \to
\mathbb{P}(V_\lambda)$ be the projective embedding determined by
this choice of basis.  (Note that one typically embeds
$G/P_\lambda$ into $\mathbb{P}(V_\lambda^\ast)$ as there is no
need for a choice of basis, but it is more convenient for us to
embed into $\mathbb{P}(V_\lambda)$.)

The following theorem has been proven by R. Dabrowski for certain
\emph{generic} $T$-orbits in $G/P$ for $G$ an arbitrary
semi-simple complex Lie group, see \cite{Dabrowski}.  Herein lies
the first proof for \emph{arbitrary} $T$-orbits in the case $G =
\SL(n,\C)$.

\begin{theorem}
The Zariski closure of any $T$-orbit in $\SL(n,\C)/P_\lambda
\hookrightarrow \mathbb{P}(V_\lambda)$ is a projectively normal
toric variety.
\end{theorem}

\begin{proof}
Let $x \in \SL(n,\C)/P_\lambda \subset \mathbb{P}(V_\lambda)$. Let
$cl(T \cdot x)$ denote the Zariski closure of the orbit $T \cdot
x$. Let $\t{Aff}(cl(T \cdot x)) \subset V_\lambda$ denote the
associated affine cone; it is easy to see that $\t{Aff}(cl(T \cdot
x)) = cl(H \cdot v_x)$ where $v_x$ is any nonzero vector on the
line $x$, since the scalar matrices in $H$ fill out all nonzero
multiples of points in $T \cdot v_x$.

Let $g \in \SL(n,C)$ be such that $gP_\lambda = x$.  Now
$wt_{\widetilde{\lambda}}(g) = A(v_x)$.  Hence by Proposition
$\ref{orbitprop}$, the affine toric variety $\t{Aff}(cl(T \cdot
x))$ is normal if and only if the semigroup algebra
$\C[\N(wt_{\widetilde{\lambda}}(g))]$ is normal, which we have
already shown.  This means that the projective toric variety $cl(T
\cdot x)$ is projectively normal.
\end{proof}

\end{document}